\documentclass{amsart}
\usepackage{amsmath}
\usepackage{amscd}
\usepackage{amssymb}
\usepackage{amsfonts}
\newtheorem{theorem}{Theorem}[section]
\newtheorem{lemma}[theorem]{Lemma}
\newtheorem{corollary}[theorem]{Corollary}

\theoremstyle{definition}

\theoremstyle{remark}
\newtheorem{remark}[theorem]{Remark}
\numberwithin{equation}{section}

\begin{document}
\title[Relating diameter and mean curvature for Riemannian
submanifolds] {Relating diameter and mean curvature for\\ Riemannian
submanifolds}

\author{Jia-Yong Wu}
\address{Department of Mathematics, Shanghai Maritime University,
Haigang Avenue 1550, Shanghai 201306, P. R. China}
\email{jywu81@yahoo.com}

\author{Yu Zheng}
\address{Department of Mathematics, East China Normal
University, Dong Chuan Road 500, Shanghai 200241, P. R. China}
\email{zhyu@math.ecnu.edu.cn}

\thanks{This work is partially supported by the NSFC10871069.}


\date{December 20, 2009.}

\dedicatory{}


\begin{abstract}
Given an $m$-dimensional closed connected Riemannian manifold $M$
smoothly isometrically immersed in an $n$-dimensional Riemannian
manifold $N$, we estimate the diameter of $M$ in terms of its mean
curvature field integral under some geometric restrictions, and
therefore generalize a recent work of Topping in the Euclidean case
(Comment. Math. Helv., 83 (2008), 539--546).
\end{abstract}
\maketitle

\section{Introduction}\label{sec1}
Let $M\rightarrow N$ be an isometric immersion of Riemannian
manifolds of dimension $m$ and $n$, respectively. In this paper, we
estimate the intrinsic diameter of the closed submanifold $M$ in
terms of its mean curvature vector integral, under some geometric
restrictions involving the volume of $M$, the sectional curvatures
of $N$ and the injectivity radius of $N$. In particular, we can
estimate the intrinsic diameter of the closed submanifold $M$ in
terms of its mean curvature vector integral without any geometric
restriction, provided the sectional curvatures of the ambient
manifold $N$ is non-positive. Our work was inspired by the following
result of P. Topping~\cite{[Topping2]} who treated the case
$N=\mathbb{R}^n$.

\vspace{0.5em}

\noindent \textbf{Theorem A.} (P. Topping~\cite{[Topping2]})
\emph{For $m\geq 1$, suppose that $M$ is an $m$-dimensional closed
(compact, no boundary) connected manifold smoothly immersed in
$\mathbb{R}^n$. Then there exists a constant $C(m)$ dependent only
on $m$ such that its intrinsic diameter $d_{int}$ and mean curvature
$H$ are related by
\begin{equation}\label{inq1}
d_{int}\leq C(m)\int_M |H|^{m-1}d\mu,
\end{equation}
where $d_{int}:=\max_{x,y\in M}dist_M(x,y)$ and $\mu$ is the measure
on $M$ induced by the ambient space. In particular, we can take
$C(2)=\frac{32}{\pi}$.}

\vspace{0.5em}

Prior to the Topping's work, L. Simon in \cite{[Simon]} (see also
\cite{[Topping0]}) derived an interesting estimate of the external
diameter $d_{ext}:=\max_{x,y\in
M^2\hookrightarrow\mathbb{R}^3}|x-y|_{\mathbb{R}^3}$ of a closed
connected surface $M^2$ immersed in $\mathbb{R}^3$ in terms of its
area and Willmore energy:
\begin{equation}\label{ext1}
d_{ext}<\frac 2\pi\left(Area(M^2)\cdot\int_{M^2}
|H|^2d\mu\right)^\frac 12.
\end{equation}
At the core of the proof of \eqref{ext1} is the following assertion
that one cannot simultaneously have small area and small mean
curvature in a ball within the surface. In other words, for all
$r>0$, we have
\begin{equation}\label{ext2}
\pi\leq\frac {A_{ext}(x,r)}{r^2}+\frac 14\int_{B_{ext}(x,r)}
|H|^2d\mu,
\end{equation}
where $B_{ext}(x,r)$ and $A_{ext}(x,r)$ denote the subset of $M^2$
immersed inside the open extrinsic ball in $\mathbb{R}^3$ centred at
$x$ of radius $r>0$ and its area, respectively. This type of
estimate is from \cite{[Simon]}, and with these sharp constants from
\cite{[Topping0]}. Combining this fact with a simple covering
argument, one can derive \eqref{ext1}. Note that if $M^2$ is a
surface of constant mean curvature $H$ immersed in $\mathbb{R}^3$,
P. Topping in ~\cite{[Toppinga]} used a different method and
established the following inequality
\[
d_{ext}\leq \frac{A|H|}{2\pi}.
\]
Equality is achieved when $M^2$ is a sphere.

Following the idea of L. Simon's proof, P. Topping in
\cite{[Topping2]} proved Theorem A by considering a refined version
of \eqref{ext2} for any dimensional manifold immersed in
$\mathbb{R}^n$. Roughly speaking, P. Topping asserted that the
maximal function and volume ratio (see their definitions in
Section~\ref{sec2}) cannot be simultaneously smaller than a fixed
dimensional constant. This assertion can be confirmed by means of
the Michael-Simon Sobolev inequality for submanifolds of Euclidean
space~\cite{[Michael-Simon]}. Then using this assertion and a
covering lemma, one can derive \eqref{inq1} immediately. As an
application, H.-Z. Li in a recent paper~\cite{[Li]} used Theorem A
to discuss the convergence of the volume-preserving mean curvature
flow in Euclidean space under some initial integral pinching
conditions.

On the other hand, as we all known, D. Hoffman and J. Spruck
in~\cite{[Hoffman-Spruck]} extended the Michael-Simon
result~\cite{[Michael-Simon]} to a general Sobolev inequality for
submanifolds of a Riemannian manifold under some geometric
restrictions. To formulate their result, we need some notations
in~\cite{[Hoffman-Spruck]}. Let $M\rightarrow N$ be an isometric
immersion of Riemannian manifolds of dimension $m$ and $n$,
respectively. We denote the sectional curvatures of $N$ by $K_N$.
The mean curvature vector field of the immersion is given by $H$. We
write $\bar{R}(M)$ for the injectivity radius of $N$ restricted to
$M$ (or minimum distance to the cut locus in $N$ for all points in
$M$). Let us denote by $\omega_m$ the volume of the unit ball in
$\mathbb{R}^m$ and let $b$ be a positive real number or a pure
imaginary one.

\vspace{0.5em}

\noindent \textbf{Theorem B.} (D. Hoffman and J.
Spruck~\cite{[Hoffman-Spruck]}) \emph{Let $M\rightarrow N$ be an
isometric immersion of Riemannian manifolds of dimension $m$ and
$n$, respectively. Some notations are adopted as above. Assume
$K_N\leq b^2$ and let $h$ be a non-negative $C^1$ function on $M$
vanishing on $\partial M$. Then
\begin{equation}\label{sobolev}
\left(\int_{M} h^{m/(m-1)}d\mu\right)^{(m-1)/m}\leq c(m)\int_{M}
\left[|\nabla h|+h|H|\right]d\mu,
\end{equation}
provided
\begin{equation}\label{condition1}
b^2(1-\beta)^{-2/m}\left(\omega_m^{-1}Vol(\mathrm{supp}
h)\right)^{2/m}\leq 1
\end{equation}
and
\begin{equation}\label{condition2}
2\rho_0\leq\bar{R}(M),
\end{equation}
where
\begin{equation*}
\rho_0=\left\{
\begin{aligned}
        b^{-1}\sin^{-1}\left[b(1-\beta)^{-1/m}\left(\omega_m^{-1}Vol(\mathrm{supp}
h)\right)^{1/m}\right]\quad\quad&\mathrm{for}\,\,\, b\,\,\, \mathrm{real},\\
(1-\beta)^{-1/m}\left(\omega_m^{-1}Vol(\mathrm{supp}h)\right)^{1/m}
\quad\quad&\mathrm{for}\,\,\, b\,\,\, \mathrm{imaginary}.
\end{aligned}
\right.
\end{equation*}
Here $\beta$ is a free parameter, $0<\beta<1$, and
\begin{equation}\label{xishu}
c(m):=c(m,\beta)=\pi\cdot
2^{m-1}\beta^{-1}(1-\beta)^{-1/m}\frac{m}{m-1}\omega_m^{-1/m}.
\end{equation}}

\begin{remark}
In Theorem B, we may replace the assumption $h\in C^1(M)$ by $h\in
W^{1,1}(M)$. As the mentioned remark in~\cite{[Hoffman-Spruck]}, the
optimal choice of $\beta$ to minimize $c$ is $\beta=m/(m+1)$. When
$b$ is real we may replace condition~\eqref{condition2} by the
stronger condition $\bar{R}\geq\pi b^{-1}$. When $b$ is a pure
imaginary number and the Riemannian manifold $N$ is simply connected
and complete, $\bar{R}(M)=+\infty$. Hence
conditions~\eqref{condition1} and~\eqref{condition2} are
automatically satisfied.
\end{remark}

Motivated by the work of P. Topping, it is natural to expect that
there exists a general geometric inequality for submanifolds of a
Riemannian manifold, which is similar to Theorem A. Fortunately,
following closely the lines of the Topping's proof of Theorem A  in
\cite{[Topping2]}, we can employ the Hoffman-Spruck Sobolev
inequality for submanifolds of a Riemannian manifold together with a
covering lemma to derive the desired results.
\begin{theorem}\label{main}
For $m\geq 1$, suppose that $M$ is an $m$-dimensional closed
connected Riemannian manifold smoothly isometrically immersed in an
$n$-dimensional complete Riemannian manifold $N$ with $K_N\leq b^2$.
For any $0<\alpha<1$, if
\begin{equation}\label{maincond1}
b^2(1-\alpha)^{-2/m}\left(\omega_m^{-1}Vol(M)\right)^{2/m}\leq 1
\end{equation}
and
\begin{equation}\label{maincond2}
2\rho_0\leq\bar{R}(M),
\end{equation}
where
\begin{equation*}
\rho_0=\left\{
\begin{aligned}
b^{-1}\sin^{-1}\left[b(1-\alpha)^{-1/m}\left(\omega_m^{-1}Vol(M)\right)^{1/m}\right]
\quad\quad&\mathrm{for}\,\,\, b\,\,\, \mathrm{real},\\
(1-\alpha)^{-1/m}\left(\omega_m^{-1}Vol(M)\right)^{1/m}
\quad\quad&\mathrm{for}\,\,\, b\,\,\, \mathrm{imaginary},
\end{aligned}
\right.
\end{equation*}
then there exists a constant $C(m,\alpha)$ dependent only on $m$ and
$\alpha$ such that
\[
d_{int}\leq C(m,\alpha)\int_M |H|^{m-1}d\mu.
\]
In particular, we can take
$C(2,\alpha)=\frac{576\pi}{\alpha^2(1-\alpha)}$.
\end{theorem}

\begin{remark}\label{rema1}
In Theorem \ref{main}, the coefficients $C(m,\alpha)$ are not
identical to (but strongly dependent on) the coefficients $c(m)$ in
Theorem B. From \eqref{budesh2} and \eqref{zuihou} we can find that
$C(m,\alpha)$ can still arrive at the minimum, when
$\alpha=\frac{m}{m+1}$. The conditions of \eqref{maincond1} and
\eqref{maincond2} are similar to the restrictions of
\eqref{condition1} and \eqref{condition2} in Theorem B, and they
guarantee that the Hoffman-spruck Sobolev for submanifolds of a
Riemannian manifold can be applied in the proof of our theorem. When
$b$ is real we may replace condition~\eqref{maincond2} by the
stronger condition $\bar{R}\geq\pi b^{-1}$. When $b$ is a pure
imaginary number and the Riemannian manifold $N$ is simply connected
and complete, $\bar{R}(M)=+\infty$, and hence
conditions~\eqref{maincond1} and~\eqref{maincond2} are automatically
satisfied.
\end{remark}

In particular, when $N=\mathbb{R}^n$, $K_N\equiv0$ and
$\bar{R}(M)=+\infty$, and hence there are also no volume
restrictions on $M$. Combining this with Remark \ref{rema1}, if $b$
is pure imaginary or zero,  then we see that
conditions~\eqref{maincond1} and~\eqref{maincond2} are automatically
satisfied, and hence we conclude that
\begin{corollary}\label{cor1}
For $m\geq 1$, suppose that $M$ is an $m$-dimensional closed
connected Riemannian manifold smoothly isometrically immersed in an
$n$-dimensional simply connected, complete, nonpositively curved
Riemannian manifold $N$ ($K_N\leq 0$). For any $0<\alpha<1$, then
there exists a constant $C(m,\alpha)$ dependent only on $m$ and
$\alpha$ such that
\[
d_{int}\leq C(m,\alpha)\int_M |H|^{m-1}d\mu,
\]
where $\min_{0<\alpha<1}C(m,\alpha)=C(m,\frac{m}{m+1})$. In
particular, we can take $\min_{0<\alpha<1}C(2,\alpha)=C(2,\frac
23)=3888\pi$.
\end{corollary}

We remark that the constants $C(2,\alpha)$ in Theorem~\ref{main} and
Corollary~\ref{cor1} are not optimal in general. The proof of
Theorem~\ref{main} follows the proof in the Euclidean
case~\cite{[Topping2]}. Theorem~\ref{main} and Corollary~\ref{cor1}
may have many interesting applications which we have not discussed
here. For example, we may borrow Li's idea of~\cite{[Li]} and apply
our Theorem~\ref{main} to study the convergence problem of the
volume-preserving mean curvature flow in Riemannian manifolds. We
will explore this aspect in the future.

Besides the above works, the closest precedent for our theorem is
another P. Topping's work on diameter estimates for intrinsic
manifolds evolving under the Ricci flow~\cite{[Topping]}. In the
Ricci flow case, P. Topping explored a log-Sobolev inequality of the
Ricci flow (see Theorem 3.4 in~\cite{[Topping]}), which can be
derived by the monotonicity of Perelman's $\mathcal{W}$-functional
(see~\cite{[CLN]},~\cite{[Perelman]},~\cite{[Topping1]}). However a
core tool of proving Theorem~\ref{main} is the Hoffman-Spruck
Sobolev inequality.

The rest of this paper is organized as follows. In
Section~\ref{sec2}, we will prove Lemma~\ref{lemm1}. The proof needs
the key Hoffman-Spruck Sobolev inequality. In Section~\ref{sec3}, we
will finish the proof of Theorem~\ref{main} using Lemma~\ref{lemm1}
of Section~\ref{sec2} and a covering lemma.

\section{Estimates for maximal function and volume ratio}
\label{sec2}

In this section we first introduce two useful geometric quantities:
the maximal function and the volume ratio. Then we apply the
Hoffman-Spruck Sobolev inequality to prove the following important
Lemma~\ref{lemm1}, which is essential in the proof of
Theorem~\ref{main}.

Given $x\in M^m$, with respect to a given metric, we denote the open
geodesic ball in $M^m$ centred at $x$ and of intrinsic radius $r>0$
by $B(x,r)$, and its volume by
\[
V(x,r):=Vol(B(x,r)).
\]
Following Topping's definitions in~\cite{[Topping2]}, when $m\geq
2$, we introduce the \emph{maximal function}
\begin{equation}\label{defin1}
M(x,R):=\sup_{r\in(0,R]}r^{-\frac{1}{m-1}}[V(x,r)]
^{-\frac{m-2}{m-1}}\int_{B(x,r)}|H|d\mu
\end{equation}
and the \emph{volume ratio}
\begin{equation}\label{defin2}
\kappa(x,R):=\inf_{r\in(0,R]}\frac{V(x,r)}{r^m}
\end{equation}
for any $R>0$.

Similar to Lemma 1.2 in~\cite{[Topping2]}, we have the following
general result.
\begin{lemma}\label{lemm1}
For $m\geq 2$, suppose that $M$ is an $m$-dimensional Riemannian
manifold smoothly isometrically immersed in an $n$-dimensional
Riemannian manifold $N$ with $K_N\leq b^2$, which is complete with
respect to the induced metric. For any $0<\alpha<1$, if
conditions~\eqref{maincond1} and~\eqref{maincond2} are satisfied,
then there exists a constant $\delta>0$ dependent only on $m$ and
$\alpha$ such that for any $x\in M$ and $R>0$, at least one of the
following is true:
\begin{enumerate}
  \item $M(x,R)\geq\delta$;
\item $\kappa(x,R)>\delta$.
\end{enumerate}
In the case of closed surfaces ($m=2$) in $N$, we can choose
$\delta=\frac{\alpha^2(1-\alpha)}{144\pi}$.
\end{lemma}

\begin{remark}
In Lemma~\ref{lemm1}, when $b$ is a pure imaginary number and the
Riemannian manifold $N$ is simply connected and complete,
$\bar{R}(M)=+\infty$. Hence conditions~\eqref{maincond1}
and~\eqref{maincond2} are automatically satisfied.
\end{remark}

Now we will finish the proof of Lemma~\ref{lemm1}.

\begin{proof}[Proof of Lemma~\ref{lemm1}]
We follow the ideas of the proof of Lemma 1.2 in~\cite{[Topping2]}.
Suppose that $M(x,R)<\delta$ for some constant $\delta>0$, which
will be chosen later. According to the definition of the maximal
function $M(x,R)$, we know that for all $r\in(0,R]$
\begin{equation}\label{equ1}
\int_{B(x,r)}|H|d\mu<\delta
r^{\frac{1}{m-1}}[V(x,r)]^{\frac{m-2}{m-1}}.
\end{equation}
Note that for fixed $x$, $V(r):=V(x,r)$ is differentiable for almost
all $r>0$. For such $r\in(0,R]$, and any $s>0$, we define a
Lipschitz cut-off function $h$ on $M$ by
\begin{equation}\label{constr}
h(y)=\left\{
\begin{aligned}
1 \quad\quad& y\in B(x,r)\\
1-\frac{1}{s}(dist_M(x,y)-r) \quad\quad& y\in B(x,r+s)\setminus B(x,r)\\
0 \quad\quad& y\not\in B(x,r+s).
\end{aligned}
\right.
\end{equation}

Since function $\sin^{-1}x$ is increasing on $[0,1]$ and
$Vol(\mathrm{supp} h)\leq Vol(M)$, we easily see that
conditions~\eqref{maincond1} and~\eqref{maincond2} guarantee the
function $h$ of~\eqref{constr} to satisfy
conditions~\eqref{condition1} and~\eqref{condition2}, where
$\beta=\alpha$. Substituting this function to the Hoffman-Spruck
Sobolev inequality from Theorem B, we derive that
\begin{equation*}
\begin{aligned}
V(r)^{(m-1)/m}&\leq\left(\int_{M}
h^{m/(m-1)}d\mu\right)^{(m-1)/m}\\
&\leq c(m)\left[\frac{V(r+s)-V(r)}{s}+\int_{B(x,r+s)}|H|d\mu\right],
\end{aligned}
\end{equation*}
where $c(m):=c(m,\alpha)=\pi\cdot
2^{m-1}\alpha^{-1}(1-\alpha)^{-1/m}\frac{m}{m-1}\omega_m^{-1/m}$.
Letting $s\downarrow 0$, we conclude that
\begin{equation*}
\begin{aligned}
V(r)^{(m-1)/m}\leq
c(m)\left[\frac{dV}{dr}+\int_{B(x,r)}|H|d\mu\right].
\end{aligned}
\end{equation*}
Combining this with~\eqref{equ1}, we have
\begin{equation}\label{budesh1}
\frac{dV}{dr}+\delta r^{\frac{1}{m-1}}V(r)^{\frac{m-2}{m-1}}-
c(m)^{-1}V(r)^{\frac{m-1}{m}}>0.
\end{equation}
Now we assume that $\delta>0$ is sufficiently small so that
$\delta<\omega_m$, and define another smooth function
\[
v(r):=\delta r^m.
\]
Then a straightforward computation yields
\begin{equation}\label{budes}
\frac{dv}{dr}+\delta r^{\frac{1}{m-1}}v(r)^{\frac{m-2}{m-1}}-
c(m)^{-1}v(r)^{\frac{m-1}{m}}=\left(m\delta+\delta^{\frac{2m-3}{m-1}}
-c(m)^{-1}\delta^{\frac{m-1}{m}}\right)r^{m-1}.
\end{equation}
We can see that
\begin{equation}\label{budesh2}
\frac{dv}{dr}+\delta r^{\frac{1}{m-1}}v(r)^{\frac{m-2}{m-1}}-
c(m)^{-1}v(r)^{\frac{m-1}{m}}\leq0
\end{equation}
as long as $\delta>0$ is sufficiently small, depending only on $m$
and $\alpha$.

Notice the fact that $V(r)/r^m\rightarrow \omega_m$ as $r\downarrow
0$, while $v(r)/r^m=\delta<\omega_m$. And combining
inequalities~\eqref{budesh1} and~\eqref{budesh2}, we conclude that
\[
V(r)>v(r)
\]
for all $r\in (0,R]$. Otherwise, there exists a fixed $r_0$ such
that $V(r_0)=v(r_0)$ and $V(r)>v(r)$ for all $r\in(0,r_0)$. Then
from~\eqref{budesh1} and~\eqref{budesh2}, we can derive
\[
\frac{dV}{dr}\Big|_{r=r_0}>\frac{dv}{dr}\Big|_{r=r_0}.
\]
Namely,
\[
\frac{dV}{dr}>\frac{dv}{dr}
\]
in any sufficiently small neighborhood of $r_0$, which is impossible
since $V(r_0)=v(r_0)$ and $V(r)>v(r)$ for all $r\in(0,r_0)$.

Therefore
\[
\kappa(x,R):=\inf_{r\in(0,R]}\frac{V(x,r)}{r^m}>\delta,
\]
which completes the proof of Lemma~\ref{lemm1}.

In the case of closed surfaces ($m=2$) in $N$, we can choose
$\delta=\frac{c(2,\alpha)^{-2}}{9}=\frac{\alpha^2(1-\alpha)}{144\pi}$
to satisfy \eqref{budesh2} and the constraint condition
$\delta<\omega_2=\pi$.
\end{proof}

\section{Diameter Control}\label{sec3}
In this section we can follow the lines of~\cite{[Topping2]}
or~\cite{[Topping]}, and easily prove Theorem~\ref{main} by using
Lemma~\ref{lemm1} and a covering lemma. For the completeness of this
paper, here we still give the detailed proof of Theorem~\ref{main}.
\begin{proof}[Proof of Theorem~\ref{main}]
We may assume $m\geq 2$ since the case $m=1$ is trivial. Now we
choose $R>0$ sufficiently large so that the total volume of the
closed manifold $M$ is less than $\delta R^m$, where $\delta$ is
given by Lemma~\ref{lemm1} (Notice that $\delta$ does not depend on
$R$). In particular, for all $z\in M$, we must have
\[
\kappa(z,R)\leq \frac{V(z,R)}{R^m}\leq \delta.
\]
Hence by Lemma~\ref{lemm1}, as long as conditions~\eqref{maincond1}
and~\eqref{maincond2} are satisfied, we must have the maximal
function $M(x,R)\geq\delta$. Namely, there exists $r=r(z)$ such that
\begin{equation}
\begin{aligned}\label{inequaty1}
\delta&\leq r^{-\frac{1}{m-1}}V(z,r)^{-\frac{m-2}{m-1}}\int_{B(z,r)}|H|d\mu\\
&\leq
r^{-\frac{1}{m-1}}\left(\int_{B(z,r)}|H|^{m-1}d\mu\right)^{\frac{1}{m-1}},
\end{aligned}
\end{equation}
where we used the H\"{o}lder inequality for the second inequality
above. Hence
\begin{equation}\label{inequaty2}
r(z)\leq\delta^{1-m}\int_{B(z,r(z))}|H|^{m-1}d\mu.
\end{equation}

Now we have to pick appropriate points $z$ at which to
apply~\eqref{inequaty2}. Let $x_1, x_2\in M$ be extremal points in
$M$. This means that $d_{int}=dist_M(x_1,x_2)$. Let $\Sigma$ be a
shortest geodesic connecting $x_1$ and $x_2$. Obviously, $\Sigma$ is
covered by the balls $\{B(z,r(z)):z\in \Sigma\}$. By a modification
of the covering lemma (see Lemma 5.2 in~\cite{[Topping]}), there
exists a countable (possibly finite) set of points
$\{z_i\in\Sigma\}$ such that the balls $\{B(z_i,r(z_i))\}$ are
disjoint, and cover at least a fraction $\rho$, where
$\rho\in(0,\frac 12)$ of $\Sigma$:
\[
\rho d_{int}\leq \sum_i 2r(z_i).
\]
Combining this with~\eqref{inequaty2}, we have
\begin{equation*}
\begin{aligned}
d_{int}&\leq\frac{2}{\rho}\sum_ir(z_i)\\
&\leq\frac{2}{\rho}\delta^{1-m}\sum_i\int_{B(z_i,r(z_i))}|H|^{m-1}d\mu\\
&\leq \frac{2}{\rho}\delta^{1-m}\int_M|H|^{m-1}d\mu,
\end{aligned}
\end{equation*}
where $\delta>0$ is sufficiently small, depending only on $m$ and
$\alpha$. Letting $\rho\to \frac 12$, we arrived at
\begin{equation}\label{zuihou}
d_{int}\leq 4\delta^{1-m}\int_M|H|^{m-1}d\mu.
\end{equation}
Hence the desired theorem follows. If $m=2$, we can choose
$4\delta^{1-m}=\frac{576\pi}{\alpha^2(1-\alpha)}$, since
$\delta=\frac{\alpha^2(1-\alpha)}{144\pi}$.
\end{proof}

\section*{Acknowledgment}
The authors thank the referee for various comments and suggestions
that helped improve this paper.


\end{document}